# Joint transfer pricing decision on tangible and intangible assets for multinational firms

Yaling Kang

School of Management, Zhejiang University of Finance & Economics, 310018, Hangzhou, Zhejiang, China

E-mail address: kangyaling622@foxmail.com

Zujun Ma

School of Management, Zhejiang University of Finance & Economics, 310018, Hangzhou, Zhejiang, China

E-mail address: zjma@zufe.edu.cn

ORCID: 0000-0002-7120-9106

Xin Tian, Corresponding Author

School of Management Science and Engineering, Dongbei University of Finance and Economics, 116025, Dalian, Liaoning, China

E-mail address: tianxin202311@163.com

Zhiqiao Wu, Corresponding Author

School of Management Science and Engineering, Dongbei University of Finance and Economics, 116025, Dalian, Liaoning, China

E-mail address: wuzhiqiao@dufe.edu.cn

ORCID: 0000-0002-7358-6175

The authors are listed in alphabetical order, the first author and the second author are listed as co-first authors in this manuscript.

# Joint transfer pricing decision on tangible and intangible assets for multinational firms

**Abstract:** While conventional multinational firms (MNFs) often avoid taxes by transferring their profits to low-tax regions through *markup* on tangible asset costs, high-tech MNFs may avoid taxes by transferring *royalty* fees to intangible assets (i.e., royalty-based transfer prices). This study investigates the effects of tax differences, markups, and royalties on decision-making. We also compare the different effects of markups and royalties on the improvement of MNFs' after-tax profit under two main business structures: the commissionaire operational structure ($C$) with complete information, and the limited-risk operational structure ($R$) in the principal–agent setting. We find that the tax difference always improves MNFs' profits under the $C$ structure, whereas non-monotonic behavior exists under the $R$ structure. More interestingly, when the order quantity is relatively small, the markup improves MNFs' profits faster than the royalty; conversely, the royalty improves MNFs' profits faster than the markup.

# 1. Introduction

## 1.1 Background and research motivation

Taxes are the largest single expense of multinational firms (MNFs) (Webber, 2011). To benefit from tax differentials worldwide, MNFs are motivated to transfer profits from high-to low-tax regions through transfer prices (Johnson 2006; Parr 2007; Huizinga & Laeven 2008; Choi et al., 2020). Previous research has shown that MNFs' subsidies result in approximately 32% lower taxes on average than similar local companies in high-tax jurisdictions (Egger et al., 2010). In 2015, more than $616 billion—nearly 40% of the profits of MNFs—was transferred to tax havens, including Ireland, Singapore, the Netherlands, the Caribbean, and Switzerland (Tørsløv et al., 2021).

Tax planning frequently utilizes transfer prices for both tangible and intangible assets. Conventionally, MNFs have employed a tax-efficient supply chain to avoid taxation through an intermediate product; for example, they might set a high cost in addition to the cost of raw materials. Jacobs Douwe Egberts (JDE), the world's largest pure-play coffee and tea group by revenue, with operations in more than 100 countries, adds a premium to the purchase price. Glaxo Canada purchased ranitidine, a patented active pharmaceutical ingredient used to combat stomach ulcers, from a related party in Switzerland, Adechsa S.A., paying more than $1,500 per kilogram from 1990 to 1993. During the same period, the market price of ranitidine was less than $300 per kilogram. Kellogg India (Singapore-based) adopting 5% cost-based transfer prices.

Compared to tangible asset profit transfer, intangible assets—especially IP—are well-suited for MNFs' tax planning due to their mobility and low migration costs (Griffith et al., 2014; Johnson, 2006; Dischinger & Riedel, 2011). A notable example is the "Double Irish arrangement", which relies on royalty-based transfer prices (Zucman, 2014; Juranek et al., 2018). This structure was central to the "Silicon Six" (Amazon, Facebook, Google, Netflix, Apple, Microsoft): their 2010–2017 effective tax rate was only 15.9% (vs. 35% headline rate), with 2012–2017 tax avoidance reaching £297 million for Facebook, £1.3 billion for Google, and £2.6 billion for Apple (Fair Tax Mark, 2019). Similar profit transfers via trademarks/brands occur at AstraZeneca, GlaxoSmithKline, and Shell (Reineke & Merten, 2021), while inter-company intangible transfers also enhance the core competence of MNFs (Johnson, 2006; Dischinger & Riedel, 2011).

More than 50% of MNFs are caught in tax disputes due to tax planning (KPMG, 2011). Tax authorities perceive MNFs of overvaluing transfer prices to shift profits, leading to frequent disputes (Dischinger & Riedel, 2011; Bauer & Langenmayr, 2013; Juranek et al., 2018; Choi et al., 2020). For tangible assets, GlaxoSmithKline settled a 14-year royalty valuation dispute with the U.S. Internal Revenue Service for $3.4 billion on September 11, 2006 (Burnett & Pulliam, 2014). Additionally, McDonald's France paid €737 million in back taxes after increasing royalties from 5% to 10% in 2009. The Double Irish arrangement is also controversial: the European Commission alleged Ireland granted Apple up to €13 billion in undue tax benefits in August 2016, but Apple won its appeal on July 15, 2020, as the Commission failed to prove a violation of the arm's length principle[1] (Fair Tax Mark, 2019; General Court of the European Union, 2020).

Both taxpayers and tax administrations face heavy administrative burdens in evaluating cross-border transactions. The core arm's length principle requires transfer prices to match unrelated firm transactions (Johnson, 2006; Huh & Park, 2013), relying on "comparability analysis"—"a comparison of a controlled transaction with an uncontrolled transaction or transactions" (OECD, 2017). However, applying this principle isn't easy: obtaining complete, explainable data from independent entities is challenging, and comparable enterprises may not exist. Thus, transfer pricing relies on reasonable estimation within an acceptable range (OECD, 2017), with royalty-based pricing even more challenging to compare than cost-based pricing (Juranek et al., 2018). Royalty-based pricing is only accepted if it falls within a reasonable profit proportion range (OECD, 2017).

MNFs often use multiple tax avoidance methods to diversify dispute risks. This study's framework examines both cost- and royalty-based transfer pricing, applicable to MNFs using either or both methods (e.g., McDonald's, Apple, Google, etc.). Excessively high markups/royalties trigger disputes, while overly low ones waste tax differential benefits—making markups and royalties critical to the operations of MNFs. However, real-world complexity adds challenges: global statutory corporate tax rates fell by 25% between 1985 and 2018, with the United States rate dropping from 35% to 21% in 2018 (Tørsløv et al., 2021); the COVID-19 pandemic also reshaped business strategies, digital adoption, and IP valuation (Deloitte, 2020). These changes raise key questions: What tax avoidance strategies boost MNFs' overall profits in tax-efficient supply chains? How should MNFs balance markups and royalties for tax planning?

Additionally, overseas branches bring tax arbitrage but require different decentralization levels (operational structures), affecting MNFs' decisions. Following the operational structure set out in Shunko et al. (2017), two common structures for profit transfer via markups and royalties are:

- **Commissionaire structure ($C$)** (full information): The procurement agent (controlled foreign corporation) has no private information; the headquarters uses mandatory contracts to ensure the expected effort level, bears corresponding effort burdens, and the agent receives a fixed salary (only responsible for raw material purchases with limited risk and decision-making authority).

- **Limited-risk structure ($R$)** (principal-agent setting): The procurement agent has decision-making power and private effort information; it bears operational risks and effort costs, receives incentive wages, and headquarters controls effort levels via incentive intensity (balancing incentive costs and income to choose optimal intensity, while the agent weighs incentive salary and effort costs to select optimal effort).

Regardless of structure, the headquarters utilizes tangible and intangible assets to achieve tax savings. For tangibles, cost-plus methods arbitrage raw materials in a low-in, high-out manner (Wu & Lu, 2018). Additionally, a certain proportion of royalty-generated profits is transferred to low-tax regions for intangibles (OECD, 2017).

## 1.2 Research questions and key findings

To illustrate the balance clearly — specifically the game equilibrium where the headquarters and overseas branches resolve conflicting goals (e.g., the headquarters pursuing tax savings vs. branches seeking reasonable incentives) under different operational structures ($C$ vs. $R$) — this study models the operational decisions related to effort level and incentive intensity. Specifically, we addressed the following research questions:

- What are the optimal decisions, including the order quantity, effort level, and intensity of incentive, for the headquarters under the $C$ and $R$ structures?
- How do the main operational parameters (i.e., tax difference, markup, and royalty) affect decision-making? Which tax avoidance method can improve the headquarters' profit?

To ensure generality, we answer these questions in a random demand environment. The sequential game of the newsvendor frame is considered for the $C$ structure, and the non-monotonic behavior occurring in the optimal decisions is analyzed. For the $R$ structure, with respect to the procurement agent's private information, the headquarters cannot directly control the procurement agent's effort. Considering stochastic demand, we establish a principal–agent framework to ensure that the procurement agent pays an expected effort level through the incentive intensity, at which point the procurement agent must bear its own effort cost. We then explore the optimal solutions and the significant effect of markups and royalties on the headquarters' after-tax profits. Interestingly, the optimal profit under the $R$ structure presents a non-monotonic tendency against the tax difference.

Our analysis yields several important findings. First, under the $C$ structure, the tendency of the effort toward the tax gap is distorted by the markup: a lower markup strengthens the effort when the tax difference increases, whereas a higher markup weakens the effort when the tax difference increases. Because of the difference in markups, the headquarters presents a different profit split proportion in the low- and high-tax regions by adjusting the effort level. When the headquarters faces a low markup, even if the tax difference increases, the benefit of tax savings is not sufficiently significant for the headquarters to choose to improve the effort level to maintain performance in the high-tax retail division. However, when the markup is large, an increasing tax gap deepens the degree of tax avoidance, making the headquarters willing to distribute more profits to the low-tax procurement region by reducing the effort level.

Second, under the $R$ structure, whether the markup is large or small, the increasing tax difference always reduces the effort level because of the decreasing incentive intensity. As the tax gap increases, the headquarters tends to save its incentive costs to avoid further taxes, thereby reducing incentive intensity. Furthermore, the procurement agent lowers effort levels.

Third, unlike the $C$ structure, MNF's profit under the R structure exhibits a non-monotonic behavior when tax differences increase. An increase in the tax gap reduces the effort level of the procurement agent and compresses the performance of the retail division. When the tax gap is relatively small, tax savings cannot cover the profit compression of the retail division. As a result, the headquarters' after-tax profits show a decreasing trend, and the tax difference increases. The opposite is true when the tax gap is relatively large.

**1.3 Contribution Statement and Article's Structure**

This study falls within the field of tax-efficient supply chains, which considers the operational structures of MNFs in relation to information symmetry and asymmetry. Our research primarily contributes to the literature in three key aspects.

First, we develop a dual tax avoidance mode, expanding traditional single-mode research. Traditional studies primarily focus on a single profit transfer mode (e.g., cost-plus markups for tangible assets, such as raw materials). In contrast, this study expands the analytical scope by proposing a dual tax avoidance mode that mirrors actual corporate practices: it incorporates both cost-based transfer price (for tangible assets, such as adjusting markups on intermedia goods like green coffee beans in supply chains) and royalty-based transfer price (for intangibles, including intellectual property, brand rights, or proprietary technologies). By integrating these two mechanisms, the framework captures the synergies between tangible asset flows and intangible value transfers—an interplay that is critical to modern MNFs' tax strategies but is underrepresented in traditional single-mode research.

Second, we link random demand, operational structures, and parameter impacts to analyze the effects on profit. Existing literature overlooks the interplay between random market demand, operational structures, and tax-related parameters in the MNF research. This study addresses this issue by analyzing how cost-plus markups (for tangibles) and royalty rates (for intangibles)—treated as key parameters—influence the decision variables of all participants (e.g., branches' effort levels and headquarters' incentive intensity) under random demand, thereby maximizing the overall profits of MNFs. Additionally, it examines how cross-region tax differentials impact these decision variables and the resulting profitability of MNFs, establishing a comprehensive link between tax disparities, operational decisions, and profit outcomes.

Third, we compare tax difference impacts across operational structures and differentiate profit improvements by markup and royalty. This study examines how cross-region tax differentials affect MNFs differently under two operational structures: $C$ structure and $R$ structure, specifying variations in their influence on core decision variables (e.g., branches' effort levels and headquarters' incentive intensity) and resultant profitability. Concurrently, it distinguishes the profit improvement effects of cost-plus markups (for tangibles) and royalty rates (for intangibles), quantifying their respective contributions to MNFs' profit enhancement to clarify which parameter yields greater gains in specific contexts.

The remainder of this paper is organized as follows. Section 2 presents the literature review and our contributions to the literature. The model description and setup are presented in Section 3. In Section 4, we analyze the optimal solutions and examine the impact of operational parameters on the decision variables. Finally, Section 5 concludes the paper.

## 2. Literature Review

This study is primarily related to two streams of literature: the core features of a tax-efficient supply chain and the information structure (symmetry or asymmetry), each of which we review below.

The structure of tax-efficient supply chains has been extensively studied (Mu et al., 2022). MNFs comprise downstream distributors (Shunko et al., 2017), retailers (Shunko & Gavirneni, 2007; Huh & Park, 2013; Kim et al., 2018; Wu & Lu, 2018; Hsu et al., 2019; Lu & Wu, 2020), upstream procurement centers (Wu & Lu, 2018), and manufacturing departments (Shunko & Gavirneni, 2007; Huh & Park, 2013; Shunko et al., 2014; Kim et al., 2018; Hsu et al., 2019). Considering the previously identified core features of a tax-efficient supply chain, the MNFs in this study make decisions based on overall profits, whereas the downstream division is an independent decision-maker. This research follows a previous study setting, based on which the impact of information structure on MNFs was investigated.

Existing studies on tax-efficient supply chains can be categorized into two scenarios: information symmetry and information asymmetry. Researchers have primarily used the framework of symmetric information. For example, Shunko and Gavirneni (2007) investigated the impact of transfer pricing on the profits of MNFs. They found that the transfer price under random demand enhanced the overall profits of MNFs more than under deterministic demand. Hsu et al. (2019) studied a multinational company with both production and retail divisions in which the retail department had a rival under deterministic demand. The authors examined the influence of the tax rate on the decision to sell to the rival. In reality, random demand is a common occurrence. Lu and Wu (2020) studied a multinational company comprising a capital-constrained retail division and a manufacturing department. They explored the effect of the tax rate on the financing strategy. The effects of blockchain technology, data sharing, and real-time payment preferences have also been studied in cross-border operations (Niu et al., 2021; Niu et al., 2022; Niu et al., 2022a). In addition, several studies have investigated tariffs as value-added taxes in the context of MNFs (Xu et al., 2018; Niu et al., 2020; Yang et al., 2021; Shi et al., 2022). However, the previously mentioned studies examined the information symmetry scenario without considering information asymmetry, that is, certain private information is not known by the MNFs' headquarters during operations. In contrast, our study has a broader scope, as it addresses both information symmetry and asymmetry under different business structures.

A few studies have also investigated tax-efficient supply chains by considering asymmetric information. For example, Jung et al. (2022) investigated the impact of the arm's length principle on the profits of MNFs when private information regarding intra-company discounts is hidden from rival retailers. In addition, the effect of different transfer pricing methods on MNFs' profits under centralized and decentralized structures in a random demand environment shows that the retail division's profits in the cost-plus method are relatively high, as the profits of the MNFs in the resale price method are relatively high (Huh & Park, 2013). Huh & Park (2013) described the asymmetry of market demand information, wherein departments closer to the market have more demand information. Shunko et al. (2014) examined the effects of different transfer pricing strategies (single or dual transfer pricing systems) on outsourcing and offshoring when MNFs face a "make or buy" problem. They found that the efficiency of tax-efficient supply chains at a single transfer price is always lower than that of a dual transfer price system. Shunko et al. (2014) considered that local operators cannot know the outsourcer's production cost, whereas we assume that the procurement agent's effort level is private information under the $R$ structure. Kim et al. (2018) also considered the asymmetry in production costs between the retail division and the outsourcing company. They investigated the influence of the arm's length principle on centralization, offshoring, and outsourcing in a global supply chain. Previous studies also considered the effort level as private information (Shunko et al., 2017; Wu & Lu, 2018). However, the studies discussed herein only considered the transfer price among departments; in this simplification, there is only one way for MNFs to avoid taxes. Specifically, our study highlights two common ways of avoiding taxation: transfer prices through markups and intellectual property fees through royalties. However, few studies have answered whether tangible assets or intellectual property are more favorable to MNFs.

## 3. Model Description

In this section, we introduce the parameters of the two operating structures (see Table 1), followed by a description of the model process. Next, we highlight the assumptions underlying our model.

Table 1 Notations.

| Abbreviations | Description |
|---|---|
| $C$ | Commissionaire operational structure |
| $R$ | Limited-risk operational structure |
| Variables | Description |
| $y$ | Order quantity of the retail division |
| $e$ | The procurement center's effort level in curtailing procurement spending |
| $b$ | Intensity of incentive for the procurement center |
| Parameters | Description |
| $m$ | Exogenous retail price |
| $\gamma_0$ | Initial price of raw material |
| $\tau$ | Retail division tax rate |
| $\tau_0$ | Procurement agent tax rate |
| $\alpha$ | Markup on cost |
| $\beta$ | Royalty for the fraction of retail profit allocated in the low-tax region |
| $\eta$ | Effect of effort |
| $k$ | Cost per unit of effort |
| $D$ | Random demand |
| $\mu$ | Mean value of demand |
| $\sigma_0^2$ | Variance of demand |
| $a$ | Reservation wage |
| Functions | Description |
| $\Delta\tau$ | Tax difference between two divisions, $\Delta\tau = \tau - \tau_0$ |
| $f(d)$ | Probability density function |
| $F(d)$ | Cumulative distribution function, let $\overline{F}(d) = 1 - F(d)$ |
| $s(y)$ | Retail division's expected sales |
| $\gamma(e)$ | Actual purchase price |
| $\pi^R$ | Retail division's profit |
| $\pi^{PC}$ | Procurement agent's profit |
| $\pi^{HQ}$ | The headquarters' profit |
| $S(\pi^R)$ | Incentive contract |

Here, we consider that an MNF owns a retail division with a high tax rate $\tau$ and a procurement division with a low tax rate $\tau_0$. This approach provides the opportunity for tax arbitrage. The headquarters, as the MNF's head office, is responsible for operational decision-making and the control of capital flows.

Before the start of the sales season, the procurement agent, employed by the procurement division, needs to exert a necessary effort $e$ to reduce the cost of raw material per unit (i.e., $\gamma(e) = \gamma_0 - \eta e$), for example, making technical improvements to obtain a lower operational cost or negotiating with a raw material supplier to curtail the purchase price. This effort comes at a convexly increasing cost $c(e) = \frac{1}{2}ke^2$. Following the markup $\alpha$, the retail division purchases raw material from the procurement division at a transfer price (i.e., $T = (1+\alpha)\gamma(e)$) based on the cost-plus method. A higher markup implies that more profit is transferred to a region with lower tax rates. To prevent the MNF from falling into tax disputes, operating parameters are required to satisfy the following equation: $\gamma(e) \leq T \leq m$. To further save taxation, the headquarters transfers a fraction of the retail division's profit as royalty, $\beta$, to the low-tax region.

When the selling season commences, the retail division sells the products to the external market at exogenous retail price $m$. Taking the competitive landscape into account, MNFs tend to maintain stable retail prices, which are less prone to significant short-term fluctuations. This is attributable to their robust cost control capabilities and supply chain management strengths. Procter & Gamble (P&G) serves as an example. Instead of directly raising prices, it adopted a product sizing strategy. Take Tide's Eco-Box as a case in point. The suggested retail price is $19.99, the same as the current 150-oz Tide press-tap. However, the Eco-Box contains 30% less water, weighs 4 pounds less, and uses 60% less plastic[2,3]. In the tech sector, the average price of an iPhone in the United States remained stable at $799 from 2021 to 2023[4]. Moreover, driven by corporate social responsibility (CSR), pharmaceutical companies like GlaxoSmithKline (GSK) have made commitments. GSK has pledged a vaccine price freeze, ensuring that the Vaccine Alliance can access the lowest price for ten years from the time of transition[5]. In conclusion, given these factors, we consider the market price to be exogenous, while market demand is subject to random fluctuations.

Finally, random market demand $d$ is realized. In addition, the *ex-ante* cumulative distribution function for demand is $F$, where $f$ and $F^{-1}$ represent the probability density function and inverse of $F$, respectively. For tractability, there is no stock-out penalty cost during the sales season, and no salvage value is recognized at the end of the sales season. Furthermore, let $g(d) = \frac{df(d)}{1-F(d)}$ denote the increasing generalized failure rate (IGFR). In this regard, the demand distribution satisfies the increasing generalized failure rate, where $g(d)$ monotonically increases in $d$. This assumption is reasonable and satisfies most distributions, such as uniform, exponential, normal, log-normal, gamma, and some Weibull distributions (see Lariviere and Porteus (2001); Lariviere (2006).

Taking into account two business structures, i.e., the $C$ and $R$ structures (see Figure 1). In the first scenario, the effort is forced by the employment contract, and the cost of effort is covered by the headquarters. The procurement agent then only receives a fixed payment (i.e., $a \geq 0$). Here, we must emphasize that effort expenses are taxed in high tax regions to reduce the tax base. In this context, excessive effort imposed by the headquarters leads to an overly high cost, while insufficient effort damages orders from the retail division. The headquarters must balance the relationship between effort expenses and revenue advantages from the company's overall performance perspective by determining a reasonable level of effort.

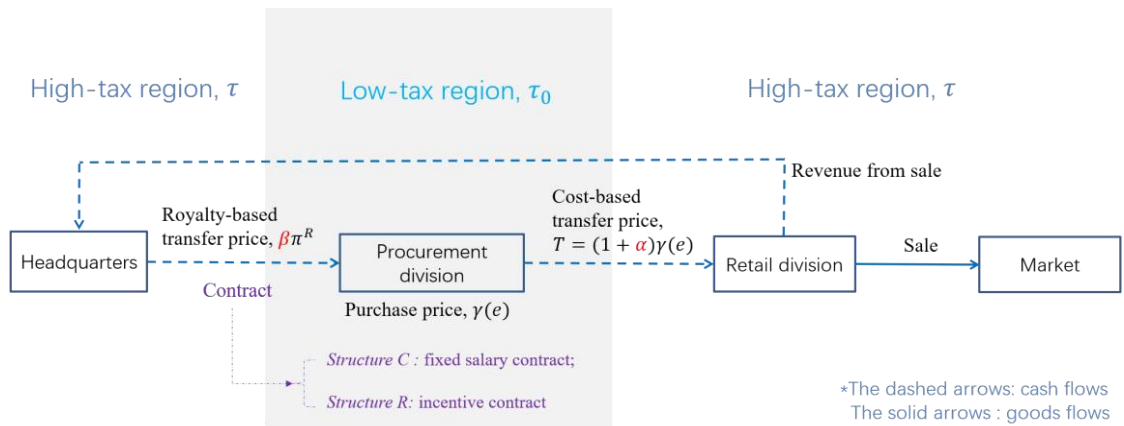

Figure 1. The MNF's business structures.

In the second scenario, moral hazard exists throughout the operation owing to asymmetric information. The procurement agent makes an effort to lower the cost of the raw material and bears the cost of their own efforts. However, accurately measuring the effort level is challenging, which can result in an overly high cost of obtaining precise information. A wise headquarters always evaluates the effort level based on the results of the effort, that is, the expected profit generated in the retail division at a certain transfer price. When facing a moral hazard, the headquarters tends not to incentivize a hard-to-observe effort level, but to pay for the outcome of the effort.

Following the standard principle–agent theory, we adopt the linear incentive contract forced on the procurement agent, which consists of a fixed payment and a bonus payment linked to the profit allocated to the procurement division: $S(\pi^R) = a + b\pi^R$, where $b$ is the intensity of the incentive set by the headquarters. At this point, the procurement agent and headquarters are two entities with different targets. Specifically, the procurement agent must weigh the trade-off between incentive wages and the costs of effort. The agent then determines the level of effort after considering its performance. As for the headquarters, there is a trade-off between the incentive cost and benefits, following which it determines the intensity of the incentive based on the company's overall performance.

Irrespective of whether it is a $C$ or $R$ structure, our setting strictly abides by the previously described tax-efficient supply chain structure, with the retail division making independent decisions. However, market demand during the sales season is uncertain. Two outcomes are possible before making an order decision: orders exceeding the demand (i.e., $d < y$) or orders that are less than the demand (i.e., $d > y$). Therefore, before the sales season, the retail division in both structures sets the inventory position decision $y$ according to its own gain $\pi^R$ after balancing the identified possibilities. Figure 2 illustrates the order of the tax-efficient supply chains.

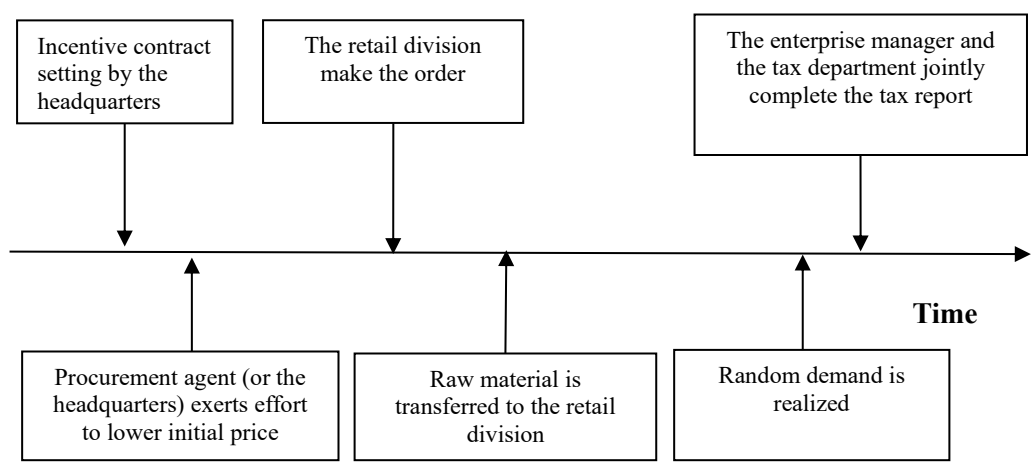

Figure 2. The sequence of the main business in a tax-efficient supply chain.

We now propose the following assumptions for simplicity:

**Assumption.** *Symmetric treatment profits and losses:* In all operational structures, we follow the arithmetic method of Shunko et al. (2017) for taxes, namely, the symmetric treatment of profits and losses.

## 4. Model Setup and Analysis

In this section, we model the $C$ and $R$ structures and analyze the effect of operational parameters (tax, markup, and royalty) on decision making.

### 4.1 Commissionaire structure ($C$)

The headquarters transfers a percentage of the retail division's profit to the procurement division and pays the procurement agent a fixed salary to exert a compulsory effort. In this scenario, the operating sequence of the MNF is as follows: First, the headquarters decides the level of effort; second, the procurement agent obtains the salary and completes the procurement function; and finally, the retail division decides the order quantity based on the transfer price. Using reverse induction, we first address the retail division problem.

$$\max_{y_C \geq 0} \pi_C^R = ms(y_C) - (1+\alpha)(\gamma_0 - \eta e_C)y_C \qquad (1)$$

where $s(y_C) = \int_0^{y_C} d\, f(d)dd + \int_{y_C}^{+\infty} y_C\, f(d)dd$—that is, the retail division's expected sales given the order quantity.

Then, according to the order reaction function from the retail division, the headquarters faces the following problem:

$$\max_{e_C \geq 0} \pi_C^{HQ} = (1-\tau)[(1-\beta)\pi_C^R - \frac{1}{2}ke_C^2] + (1-\tau_0)[\beta\pi_C^R - a + \alpha(\gamma_0 - \eta e_C)y_C] \qquad (2)$$

From Equation (2), the headquarters' after-tax profit is mainly composed of two parts: profits from the high-tax retail division and profits from the low-tax procurement division. The headquarters must allocate profits between identified parts by determining the effort level. Without loss of generality, the fixed salary was assumed to be zero. The effort under the $C$ structure is equivalent to the headquarters' efforts due to the compulsory contract. Thus, we propose the following lemma.

***Lemma 1:*** *Under the* ***C*** *structure, the optimal decisions of each participant are*

$$\begin{cases} y_C^* = F^{-1}\left(\frac{m-(1+\alpha)\gamma(e_C^*)}{m}\right) \\ e_C^* = \frac{N\eta y_C^* + (1-\tau_0)\alpha\gamma_0 \frac{dy_C}{de_C}}{(1-\tau_0)\left(\alpha\eta\frac{dy_C}{de_C}+k\right)-(\tau-\tau_0)k} \end{cases},$$

*where* $N = [(1-\tau)(1-\beta) + (1-\tau_0)\beta](1+\alpha) - (1-\tau_0)\alpha$, *and* $\frac{dy_C}{de_C} = \frac{(1+\alpha)\eta}{mf(y_C)}$.

Lemma 1 shows the optimal responses of the retail division and the headquarters. According to Lemma 1, the transfer price decreases the order quantity and effort level, whereas the royalty transferred to the procurement division increases both the order quantity and the effort level. A higher royalty (i.e., a larger $\beta$) makes the headquarters obtain more after-tax revenue. Thus, the headquarters expects the procurement agent to exert greater effort.

***Proposition 1:*** *Under the* ***C*** *structure, when the markup is relatively low (i.e., $\alpha \leq \hat{\alpha} = \frac{\beta}{1-\beta}$), a significant tax difference will strengthen the effort; when the markup is relatively high, a significant tax difference may weaken the effort.*

Recall that the headquarters' profits come from its retail and procurement divisions and that the headquarters must weigh the proportion of profit fragments by adjusting its effort level. From Proposition 1, a relatively low markup significantly weakens the tax-saving advantage developed through the procurement division, despite the large tax difference. In this situation, the retail division's performance is crucial. Within this range, increased effort can lead to higher revenue through increased order volume. Thus, the headquarters expects the procurement agent to increase the degree of effort, although the cost of the effort also increases accordingly. We believe this scenario warrants further investigation. Nevertheless, a relatively high markup may significantly improve tax arbitrage. Within this range, a significant tax difference further enhances tax savings in the procurement division;

accordingly, the headquarters expects to allocate less profit to the retail division by reducing the effort level.

This phenomenon can be explained by Lemma 1 and Proposition 1. According to Proposition 1, an increase in the markup leads to a decrease in performance in the retail division due to the excessive transfer price. The promotion effort cannot stimulate the order quantity quite well, but leads to an expensive effort cost. Therefore, the effort level decreases accordingly. Thus, from Proposition 1, we find that a significant tax difference intensifies the responses to the markup.

To illustrate Proposition 1, Figure 4 presents a numerical example.

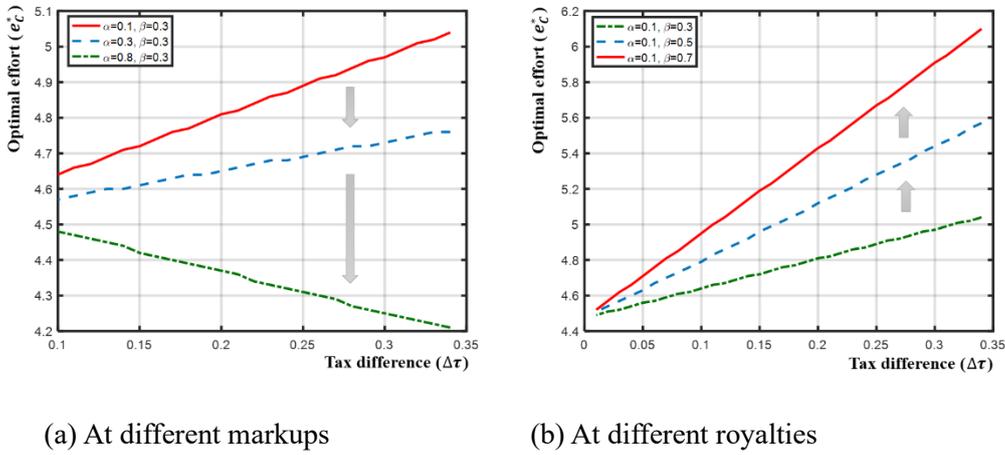

(a) At different markups       (b) At different royalties

Figure 4. The relationship between tax differences and optimal effort, where $m = 100$, $\gamma_0 = 20$, $\eta = 1$, $k = 56$, $w = 5100$, $\tau = 0.35$, $\mu = 220$, and $\sigma_0 = 30$ (The parameters are used throughout the $C$ structure).

From Figure 4(a), in the case of the thin markup ($\alpha = 0.1$), the advantage of tax arbitrage is dominated by profits from the retail division. The headquarters tends to improve order quantity by strengthening its efforts as the tax gap increases. Therefore, when the tax gap increases from $\Delta\tau = 0.05$ to $\Delta\tau = 0.3$, the headquarters' effort increases from $e_C^* = 4.56$ to $e_C^* = 4.96$. However, in the case of the high markup ($\alpha = 0.8$), the advantage of tax arbitrage becomes extremely important as the headquarters expects to weaken its efforts due to the effort cost burden. Thus, when the tax difference changes from $\Delta\tau = 0.05$ to $\Delta\tau = 0.3$, the headquarters' effort decreases from $e_C^* = 4.54$ to $e_C^* = 4.25$. Ultimately, in a given tax difference, an increasing markup (i.e., $\alpha = 0.1$, $\alpha = 0.3$, $\alpha = 0.8$) continues exacerbating the retail division's operating cost, which leads to a lowered order volume. The headquarters anticipates this change and considers a lower effort level. As Figure 4(b) indicates, the royalty proportion does not directly affect the operating costs of the retail division, unlike the markup. Therefore, an increasing royalty-based transfer price ($\beta = 0.3$, $\beta = 0.5$, $\beta = 0.7$) will encourage the headquarters to expend more effort to increase the performance of the retail division. Thus, further taxation is avoided through royalties.

***Corollary 1:*** *Under the **C** structure, the headquarters' profits have the following characteristics: the greater the tax difference, the higher the markup, the higher the royalty proportion, and the more the headquarters' profits can be improved.*

In this operating channel, the headquarters gains greater control by adjusting its effort level, enabling controlled decision-making that improves profits. An increase in the tax gap, markup, and royalty increases the convenience for headquarters in terms of tax savings, as shown in Figure 5.

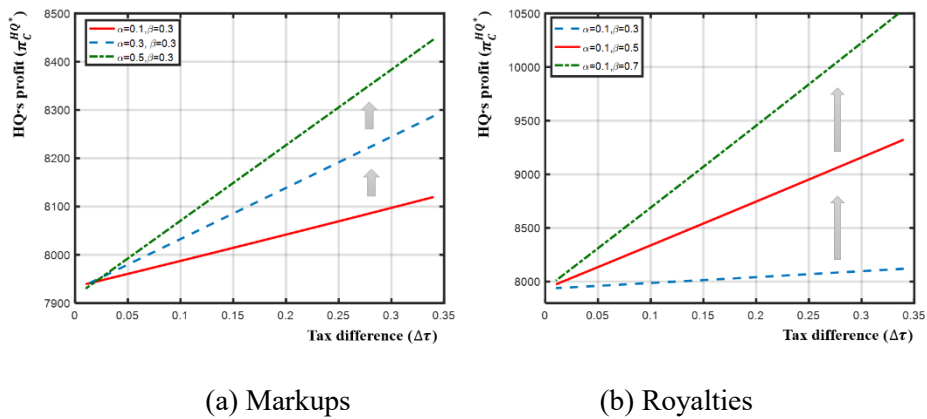

(a) Markups  (b) Royalties

Figure 5. The relationship between tax differences and the headquarters' profits.

As indicated in Figures 5(a) and (b), the headquarters can consistently generate profits from the tax gap, as the increasing tax difference continually enhances the headquarters' profits. In Figure 5(a), in the case of the same royalty ($\beta = 0.3$), with the increase of the markup ($\alpha = 0.1$, $\alpha = 0.3$, $\alpha = 0.5$), the profits transferred to the low-tax region also increase. As a result, the headquarters' profits further improve. Similarly, in Figure 5(b), in the case of the identical markup ($\alpha = 0.1$), the headquarters' transferrable profits increase as a proportion of the royalty ($\beta = 0.3$, $\beta = 0.5$, $\beta = 0.7$). In short, under the **C** structure, the tax difference, markups, and royalties can significantly enhance the headquarters' profits.

### 4.2 Limited-risk operational structure (*R*)

Under the **R** structure, the headquarters transfers a fragment of the retail division's profit to the procurement division and pays the procurement agent an incentive wage to promote further effort. In this scenario, the operating sequence of the MNF is as follows: the headquarters determines the bonus level, the procurement agent determines its own effort based on the incentive, and the retail division fulfills the order from the procurement division based on the transfer price.

Similarly, we addressed this problem using backward induction. Initially, the retail division faced a newsvendor problem.

$$\max_{y_{RN} \geq 0} \pi_{RN}^R = ms(y_{RN}) - (1 + \alpha)(\gamma_0 - \eta e_{RN})y_{RN}, \quad (3)$$

where $s(y_{RN}) = \int_0^{y_{RN}} d f(d) dd + \int_{y_{RN}}^{+\infty} y_{RN} f(d) dd$ —that is, the retail division's expected sales given the order quantity.

The headquarters then formulates a linear contract based on the procurement agent's behavior, which can not only stimulate the agent's effort but also maximize the headquarters' own performance. This model can be expressed as follows:

$$\max_{0 \leq b_{RN} < 1} \pi_{RN}^{HQ} \quad (4)$$

$$s.t. \ IR \ \pi_{RN}^{PC} \geq a$$

$$IC \ e_{RN} = \underset{e_{RN} \geq 0}{argmax} \ \pi_{RN}^{PC}$$

where $\pi_{RN}^{HQ} = (1 - \tau)(1 - \beta)\pi_{RN}^R + (1 - \tau_0)[\beta \pi_{RN}^R - a - b_{RN}\pi_{RN}^R + \alpha(\gamma_0 - \eta e_{RN})y_{RN}]$, $\pi_{RN}^{PC} = a + b_{RN}\pi_{RN}^R - \frac{1}{2}k e_{RN}^2$.

As a result, we posit the following lemma:

***Lemma 2:*** *Under the **R** structure, the optimal decisions of each participant are*

$$\begin{cases} y^*_{RN} = F^{-1}\left(\dfrac{m - (1+\alpha)\gamma(e^*_{RN})}{m}\right) \\ e^*_{RN} = \dfrac{b^*_{RN}(1+\alpha)\eta y^*_{RN}}{k} \\ b^*_{RN} = \dfrac{Ng(y^*_{RN}) + (1-\tau_0)\alpha}{(1-\tau_0)(1+\alpha)g(y^*_{RN})} \end{cases}$$

Lemma 2 addresses the optimal behavior of each participant. According to Lemma 2, a higher markup lowers the order quantity and further squeezes profits in the retail division. Even with the same incentive intensity, the available incentive salary decreases, which directly reduces the procurement agent's effort level. Intelligent leadership also recognizes that increasing revenue through incentive wages cannot fully offset the associated costs of incentives. Thus, a higher markup leads to lower incentive intensity. Moreover, higher royalties generate greater tax savings in regions with lower tax rates. In this case, the headquarters will pay more attention to procurement profits and lower the incentive intensity.

**Proposition 2:** *Under the **R** structure, optimal behavior has the following characteristics: both bonus incentive $b^*_{RN}$ and effort $e^*_{RN}$ decrease in the tax difference $\Delta\tau$.*

Under the **R** structure, the headquarters transfers a fragment of the retail division's profit to the low-tax region, and the remaining part is subject to tax arbitrage after the incentive cost is paid. As the tax difference gap increases, the tax-saving advantage of the procurement division increases. Thus, the incentive intensity of the headquarters' decision-making has a negative correlation with the tax difference. Regarding the procurement agent, when the incentive wage decreases due to the tax difference, the procurement agent expects its marginal revenue to decrease and accordingly reduces its own effort level. To illustrate Lemmas 2 and Proposition 2 better, we provide a numerical example in Figure 6.

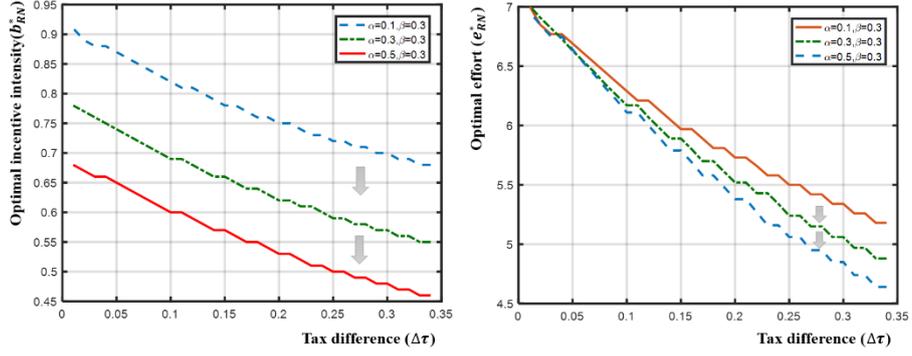

(a) At different markups

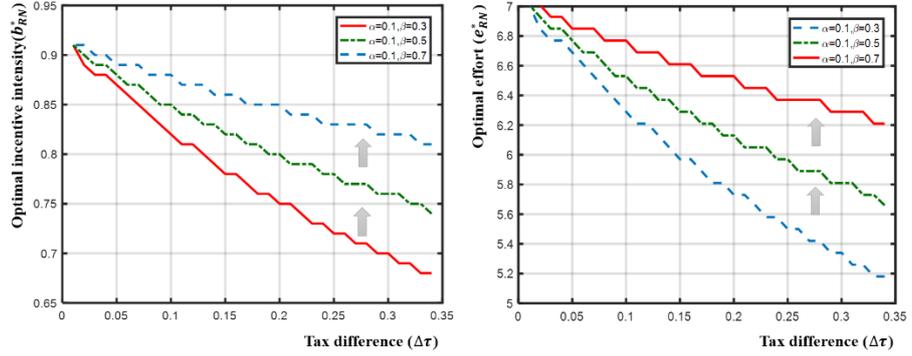

(b) At different royalties

Figure 6. Optimal decisions when changing tax differences occur, where $m = 100$, $\gamma_0 = 20$, $\eta = 1$, $k = 36$, $w = 5100$, $\tau = 0.35$, $\mu = 220$, and $\sigma_0 = 30$ (The parameters are used throughout the $C$ structure).

Based on Figure 6, in the case of $\alpha = 0.1$ and $\beta = 0.3$, as the tax gap increases, the headquarters will reduce incentive wages by lowering the incentive intensity to utilize the tax gap advantage better. Therefore, the tax gap increases from $\Delta\tau = 0.05$ to $\Delta\tau = 0.3$, and the incentive intensity decreases from $b_{RN}^* = 0.87$ to $b_{RN}^* = 0.7$. Correspondingly, the procurement agent reduces its effort from $e_{RN}^* = 6.7$ to $e_{RN}^* = 5.4$. Additionally, Figure 6(a) shows that an increase in markup leads to a decrease in incentive intensity and effort level. Taking $\Delta\tau = 0.25$ as an example, the markup from $\alpha = 0.1$ to $\alpha = 0.3$ and then to $\alpha = 0.5$ directly reduces the order quantity developed from the retail division and the margin revenue paid to the procurement agent. The procurement agent will then reduce its own effort level; that is, the effort decreases from $e_{RN}^* = 5.5$ to $e_{RN}^* = 5.3$ and then to $e_{RN}^* = 5.1$. The headquarters also realizes that a high incentive cost cannot sufficiently motivate the procurement agent, and tends to lower the incentive intensity from $b_{RN}^* = 0.73$ to $b_{RN}^* = 0.59$ and then to $b_{RN}^* = 0.5$. However, as Figure 6(b) demonstrates, increased royalties enhance the incentive intensity and effort level. Taking $\Delta\tau = 0.25$ as an example, the royalty from $\beta = 0.3$ to $\beta = 0.5$ and then to $\beta = 0.7$ promotes the profits transferred to the low-tax region. The headquarters will then pay higher incentive wages to increase incentive intensity, stimulating the procurement agent's effort and the retail division's orders. At this point, the benefits of increased incentive wages for the procurement agent outweigh the burden of the incentive wage.

***Proposition 3:*** *Under the **R** structure, the headquarters' profits have the following characteristics: (i) there exists a threshold $\Delta\tau_{RN}^{\#}$ such that if $0 < \Delta\tau < \Delta\tau_{RN}^{\#} = \tau - \tau_0^{\#}$, the headquarters' profits decrease in the tax difference $\Delta\tau$; otherwise, the headquarters' profits increase in the tax difference $\Delta\tau$ when $\Delta\tau > \Delta\tau_{RN}^{\#}$; and (ii) the higher the markup and the more the proportion of royalty, the more that the headquarters' profits can be improved.*

Recall that headquarters' performance is affected by three factors: 1) after-tax profit from the retail division, 2) tax savings on the transfer price in the low-tax region, and 3) tax savings on the royalty after the payment of the incentive wage. Based on Proposition 3, both the incentive intensity and effort level show a downward trend as the tax difference increases; the performance of the retail division also subsequently decreases in this range. When the headquarters faces a relatively low tax difference (i.e., $\Delta\tau < \Delta\tau_{RN}^{\#}$), the small tax savings cannot cover the decline in the retail division due to the increase in the tax difference. Thus, the headquarters' after-tax profits also show a downward trend.

On the other hand, when the headquarters faces a relatively high tax difference (i.e., $\Delta\tau > \Delta\tau_{RN}^{\#}$), the tax-saving advantage becomes more significant. Within this range, tax arbitrage can offset the side effects of the retail division's reduced revenue, and the headquarters' after-tax profit shows an upward trend. Therefore, the headquarters can only benefit from tax arbitrage if the tax gap is sufficiently large.

In addition, similar to the ***C*** structure, markups and royalties can always significantly promote headquarters' profits. Considering Proposition 3 and Figure 6, although an increased markup reduces the income produced by the retail division, the reduction in incentive intensity also reduces incentive costs. Thus, the headquarters' profits also increase because of an increase in markups.

On the other hand, the effect of increased royalties is transmitted through changes in the incentive intensity decision, which cannot directly affect the retail division's purchase cost, such as the markup. In other words, tax avoidance through royalties may lead to smaller chain-wide decision responses. Therefore, an increase in royalty further motivates the headquarters to provide incentives, driving the procurement agent to exert more effort to improve the retail division's profit. Finally, the headquarters' profit increases further because of higher royalties. Figure 7 illustrates Proposition 3.

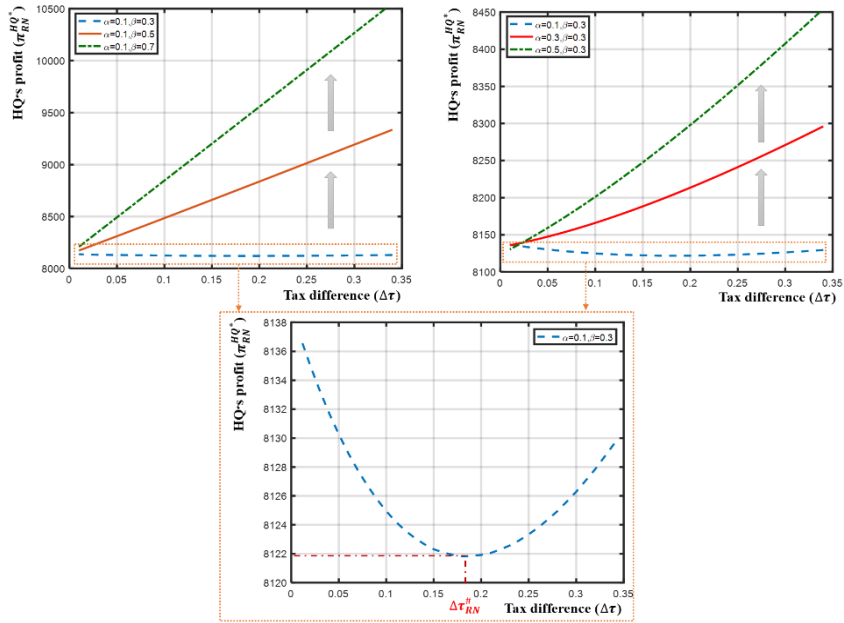

Figure 7. The relationship between the tax difference and the headquarters' profit.

From Figure 7, when $\alpha = 0.1$ and $\beta = 0.3$, there exists a threshold $\Delta\tau_{RN}^{\#} = 0.18$, such that the headquarters' profit decreases with the tax difference when the tax difference is less than the threshold (i.e., $\Delta\tau < 0.18$). When the tax difference exceeds this threshold (i.e., $\Delta\tau > 0.18$), the headquarters' profit increases with the tax gap. Considering Figure 7, we can intuitively observe that markups and royalties significantly improve the headquarters' profits. As stated in Proposition 3, when the tax difference is 0.25, the markup increases from $\alpha = 0.1$ to $\alpha = 0.3$ and then to $\alpha = 0.5$; the headquarters' profit increases from $\pi_{RN}^{HQ^*} = 8123$ to $\pi_{RN}^{HQ^*} = 8241$ and then to $\pi_{RN}^{HQ^*} = 8351.5$. In contrast, the royalty improves the headquarters' profit faster than the markup, increasing from $\beta = 0.3$ to $\beta = 0.5$ and then to $\beta = 0.7$. The headquarters' profit increases from $\pi_{RN}^{HQ*} = 8123$ to $\pi_{RN}^{HQ*} = 9014.6$ and then to $\pi_{RN}^{HQ*} = 9912$.

### 4.3 Comparison of different structures

The previous section analyzed the operational decisions under different operational structures. This section provides numerical illustrations of the profit improvement that an MNF can achieve by adopting a tax-efficient operating structure. Proposition 4 addresses the question of which choice is better for the MNF: markups or royalties.

***Proposition 4:*** *Under the **C** and **R** structures, respectively, a threshold exists, showing that the markup can improve the headquarters' profit more when the order quantity is relatively low; otherwise, the royalty can improve the headquarters' profit more when the order quantity is relatively large.*

Recall the previously stated arm's length principle dealing with royalties: the transferable profits generated through intellectual property rights account for a certain percentage of the retail division's profits. In both business structures, when the order quantity is relatively small, the meager profit developed from the retail division makes it less advantageous to save tax through royalty; thus, the markup can further improve the headquarters' profit. However, when the order quantity is relatively large, the headquarters increases the proportion of royalties to effectively avoid taxes. Specifically, markups can more easily obtain advantages under the **R** structure, whereas royalties can more easily obtain advantages under the **C** structure. This is because the order quantity under the **C** structure is always larger than that under the **R** structure, even though they have the same threshold. Thus, under the **C** structure, it is recommended to use the royalty-based transfer price method for tax planning, whereas under the **R** structure, the MNF can more reasonably avoid taxes through a cost-based transfer price. We illustrate Proposition 4 through Figure 8.

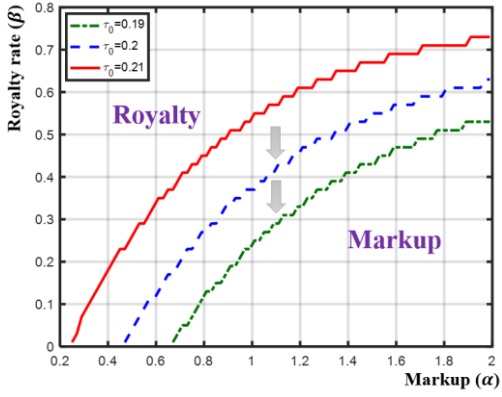 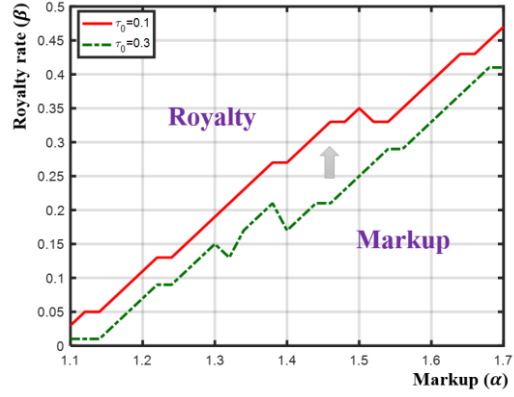

(a) $C$ structure  (b) $R$ structure

Figure 8. Profit comparison w.r.t. markup and royalty, where subfigure (a) uses parameter values $m = 100$, $\eta = 1.3$, and $k = 56$, while subfigure (b) uses $m = 100$, $\eta = 1.3$, and $k = 56$.

The lines in Figure 8 show $\frac{dE[\pi_C^{HQ^*}]}{d\alpha} = \frac{dE[\pi_C^{HQ^*}]}{d\beta}$ in different markup and royalty scenarios. As shown in the figure, the area below the curve indicates that the markup promotes headquarters' profits more efficiently, while the area above the curve shows that royalties improve headquarters' profits more efficiently. Under the $C$ structure, when the headquarters faces a set of joint transfer prices that reduce the order volume in the retail division (i.e., $g(y_C^*) < g(\hat{y}_C)$), it adopts a higher markup that avoids taxation to improve the MNF's profits effectively. However, when the headquarters face a set of joint transfer prices that increase the order quantity, a higher royalty is suggested. As shown in Figure 8(b), a similar conclusion was reached for the $R$ structure.

Figure 8(a) also confirms some of the findings in Proposition 1, where an increase in the tax gap leads to higher effort and order volume. When the low tax rate $\tau_0$ decreases from 0.21 to 0.2 and then to 0.19 (i.e., the tax difference increases), the tax advantage is not significant enough for the headquarters to increase the order quantity by enhancing efforts. Thus, royalty is more likely to prevail, and its dominant area increases. Similarly, Figure 8(b) verifies the conclusions in Proposition 2, namely that a higher tax gap reduces both the effort level and the order quantity. When the low tax rate $\tau_0$ decreases from 0.3 to 0.1 (i.e., the tax gap increases), the policy of further tax savings by lowering the incentive cost results in a lower effort level and order volume. Thus, markups are more likely to prevail, and their area of domination increases.

In addition, Figure 8 shows that both markups and royalties can continuously improve headquarters' profits in both business structures. To improve the headquarters' profits using a given tax difference ($\tau_0 = 0.2$), the markup $\alpha$ increases from 0.6 to 0.8, compared to the increase of the royalty rate $\beta$ from 0.12 to 0.28, confirming Corollary 1 and Proposition 3.

## 5. Conclusion
### 5.1 Concluding remarks and main insights

We focus on the operations of MNFs, examining tax avoidance through comparisons of markups with royalties in commissionaire and limited-risk operational structures. In a newsvendor environment, the retail division in any business structure faces inventory issues ahead of the selling season. In the commissionaire's operational structure, the headquarters adopted a compulsory contract for a procurement agent. We show that increasing the tax gap strengthens the headquarters' effort levels when the markup is relatively low. However, when the markup is high, the increased tax difference weakens the headquarters' effort levels. In a limited-risk operational structure, the headquarters applies an incentive wage to a procurement agent. Unlike the commissionaire's operational structure, the tax difference does not always improve the headquarters' after-tax profit in this case. Although the markup and royalty can continuously help the headquarters save on taxation and promote its after-tax profit in both operational structures, a threshold of order quantity exists that motivates the headquarters to adopt different preferences for markup and royalty. Furthermore, when the order quantity is relatively small (large), the markup (royalty) improves MNFs' profits much faster than the royalty (markup).

Based on our findings, we can offer the following valuable managerial insights.

- MNFs should make adjustments to decisions in a timely manner in a changing operational environment. Under a commissionaire's operational structure with tax policy adjustments resulting in a larger tax difference, it is suggested that MNFs strengthen (or weaken) their effort level only when the cost-based transfer price (i.e., markup) is relatively low (or high). Under a limited-risk operational structure with a higher tax gap, MNFs should decrease their effort levels regardless of the markup value.
- The tax gap does not always result in high after-tax profits for MNFs. The tax difference

promotes net profits only when the tax difference is relatively large.

- When a retail division produces a high volume due to a potential market, our analysis suggests that MNFs should enhance the royalty-based transfer price to improve profits more effectively. When a retail division faces fierce competition resulting in lower order quantity, our conclusions suggest that MNFs should increase the cost-based transfer price.

## 5.2 Future research

Our work has several limitations that are worth mentioning and warrant further study. First, we consider markup as a strategy parameter while overlooking its role as a decision variable. It would be interesting to examine how taxes affect strategic markup decisions. Second, cooperation between the headquarters and the procurement agent is a long-term endeavor. Thus, it is crucial to discuss how a reasonable incentive mechanism is established within the limited-risk operational structure and how long-term profits can be sustained in a repeated game.

**Notes**

[1] Article 107(1) TFEU: "Save as otherwise provided in the Treaties, any aid granted by a Member State or through State resources in any form whatsoever which distorts or threatens to distort competition by favoring certain undertakings or the production of certain goods shall, in so far as it affects trade between Member States, be incompatible with the internal market" (General Court of the European Union, 2020).

[2] https://www.nairaland.com/4835618/p-g-introduces-new-concentrated

[3] https://www.businesswire.com/news/home/20181109005291/en/PG-Thinks-Inside-the-Box-with-New-Tide-Eco-Box

[4] https://moneynotmoney.com/historical-price-of-iphone-in-united-states/

[5] https://www.gsk.com/en-gb/media/press-releases/gsk-announces-new-commitment-to-improve-access-to-vaccines-with-5-year-price-freeze-for-countries-graduating-from-gavi-alliance-support/